\documentclass[12pt]{article}
\usepackage{graphicx}
\usepackage{amsmath}
\usepackage{amsthm}
\usepackage{amssymb}
\usepackage{amsfonts}
\usepackage{epsfig}
\usepackage{epstopdf}

\def\arg {\mathop{\rm arg}\nolimits}

\def\Res{\mathop{\rm Res}}

\newtheorem{thm}{Theorem}

\numberwithin{equation}{section}

\title{{Asymptotics of  recurrence coefficients for the Laguerre weight with a singularity at the  edge}}

\author{Wu Xiao-Bo\footnote{
 {\it{E-mail
address:}} {wuxiaobo201207@163.com} }
\\
\hbox{\small \emph{Department of Mathematics, Shangrao Normal University, Shangrao, China}}
}

\date{}
\begin{document}

\maketitle


\noindent {\bf{Abstract }} In this paper, We study the asymptotics of the leading coefficients and the recurrence coefficients for the orthogonal polynomials with repect to the Laguerre weight with  singularity of root type and jump type at the soft edge via the Deift-Zhou steepest descent method. The asymptotic formulas of the leading coefficients and the recurrence coefficients for large n are described in terms of a class of analytic solutions to the the  $\sigma$-form of the Painlev\'{e} II equation and the Painlev\'{e} XXXIV equation.

\vspace{.2in} \noindent  {Keywords: Asymptotic;  Leading coefficients; Recurrence coefficients; Riemann-Hilbert approach}



\section{Introduction and statement of results} \indent\setcounter{section} {1}

\noindent

We consider the orthogonal polynomials $p_n(x)$ satisfied the following relationship
\begin{equation}\label{pLUG}
\int_0^{+\infty} p_n(x)p_m(x)w(x)dx=\delta_{n,m},\end{equation}
with the  weight perturbed by a Fisher-Hartwig singularity
\begin{equation}\label{weight}
w(x):=w(x;\alpha,\beta,\omega,\mu)=x^\alpha e^{-x}|x-\mu|^{2\beta}\left\{\begin{array}{cc}
                                                         1, & 0<x\leq\mu\\
                                                         \omega, & x>\mu
                                                       \end{array},
\right.  \end{equation}
where $\alpha>-1$, $\beta>-\tfrac 12$,  $ \omega\in \mathbb{C}\setminus(-\infty,0)$ and $\mu=4n+4^{2/3}n^{1/3}s$. It is seen that the weight $w(x)$ given as \eqref{weight} has Fisher-Hartwig singularily of root type and jump type at $x=\mu$, where $2^{2/3}n^{2/3}\left(\tfrac{\mu}{4n}-1\right)\rightarrow s$ with finite $s$ as $n\rightarrow \infty$.

Let $\gamma_n$ be the leading coefficient of the orthonormal polynomials with respect to \eqref{weight} and
\begin{equation}\label{gamma}
p_n(x)=\gamma_n\pi_n(x),
\end{equation}
then we have
the three-term recurrence relations
\begin{equation}\label{recurrence relations}
z\pi_n(x)=\pi_{n+1}(x)+a_n \pi_n(x)+b_n^2 \pi_{n-1}(x)
\end{equation}
with $\pi_{-1}(x)=0$ and $\pi_{0}(x)=1$.

When $\beta=0$ and $\omega=1$ in weight \eqref{weight}, the \eqref{weight} is the Laguerre weight; see \cite{Szego1975,Ismail}. Based on Deift-Zhou steepest descent method, the Laguerre-type weight $w(x)=x^\alpha e^{-Q(x)}$ were considered in \cite{Vanlessen}, where $Q(x)$ is a polynomial with positive leading coefficient,
and  global asymptotic expansions of the Laguerre polynomials were  studied in \cite{QiuWong}. The recurrence coefficients of orthogonal polynomials may be related to the solutions of the Painlev\'{e} equations; see \cite{ChenIts} for the Laguerre-type weight perturbed by a factor $e^{-s/x}$ with $s>0$, and  \cite{FilipukVanAsscheZhang} for semi-classical Laguerre polynomials.

The discontinuous weights were considered by the ladder operator approach in \cite{ChenGunnar}, see also \cite{BasorChen,LyuChen}  for the discontinuous Laguerre weight with a jump, and \cite{ChenFeigin,MinChen} for the Hermite weight with an isolated zero and the Gaussian weight with two jumps respectively. In \cite{LyuChen}, Lyu and Chen studied the  weight  given in \eqref{weight} in the case $\beta=0$, $\omega=0$ and $\alpha=O(n)$ as $n\rightarrow\infty$ or $\alpha$ is finite.
The logarithmic derivative of the limit of the probability of the largest eigenvalue for large $n$ was described respecting certain solutions of the Jimbo-Miwa-Okamoto $\sigma$- form of the Painlev\'{e} II equation
\begin{equation}\label{sigma form}
(\sigma'')^2+4(\sigma')^3-4s(\sigma')^2+4\sigma'\sigma- (2\beta)^2=0,
 \end{equation}
with $\beta=0$, and the asymptotics of the recurrence coefficients were also obtained.

In \cite{ItsKuijlaarsOstensson2008}, Its, Kuijlaars and \"{O}stensson investigated perturbed Gaussian unitary ensemble(pGUE) with singularity of root type at $\sqrt{2n}+\frac s{\sqrt{2}n^{1/6}}$ for bounded $s$ based on the Deift-Zhou steepest descent method. They obtained the asymptotic formulas of the limiting eigenvalue correlation kernel as $n\rightarrow\infty$ characterized in terms of a solution to the Painlev\'{e} XXXIV($\mathrm{P}_{34}$ for short) equation
\begin{equation}\label{pXXXV}
 {u}_{ss}=\frac{{u}_s^2}{2{u}}+4{u}^2+2s{u} -\frac{(2\beta)^2}{2{u}}.
  \end{equation}
Later, the pGUE with a jump at $\sqrt{2n}+\frac s{\sqrt{2}n^{1/6}}$ was studied in \cite{XuZhao} and \cite{BogatskiyClaeysIts},
and with singularity of root type and jump type at $\sqrt{2n}+\frac s{\sqrt{2}n^{1/6}}$
were investigated in \cite{WuXuZhao}. In \cite{WuXuZhao}, they found that the asymptotic formulas of the recurrence coefficients and Hankel determinant associated with the weight can be characterized in terms of a class of solutions to the equations \eqref{sigma form} and \eqref{pXXXV}. The pGUE with singularity of root type and jump type at $\sqrt{2n}+\frac s{\sqrt{2}n^{1/6}}$ for the case $2\beta\in \mathbb{N}$ and $\omega=0,1$ was  considered in \cite{ForresterWitte} earlier by applying the Okamoto $\tau$-function theory.

In this paper, we investigate the asymptotic formulas of the leading coefficients and recurrence coefficients for the perturbed Laguerre weight with singularity of root type and jump type at $4n+4^{2/3}n^{1/3}s$ defined as  \eqref{weight} by Deift-Zhou steepest decent method.
\subsection{Statement of main results}
The asymptotic formulas of leading coefficients and recurrence coefficients are characterized as follows.
\begin{thm}\label{Asymptotic of leading coefficients}
Let $\alpha>-1$, $\beta>-\frac 12$, $\omega\in \mathbb{C}\setminus (-\infty,0)$, $\mu=4n+4n^{1/3}s$, $\gamma_{n}$ be the leading coefficients in \eqref{gamma}, and $a_n$ and $b_n$ be the recurrence coefficients in \eqref{recurrence relations}. Then the following holds
\begin{gather}
\gamma_{n-1}=\frac{n^{-n-\frac \alpha 2-\beta+\frac 12}e^n}{\sqrt{2\pi}}
\left[1 +\frac{\sigma(s)}{2^{1/3}n^{1/3}} +\frac{\sigma^2(s)+2(\alpha+2\beta-1)u(s) }{2^{5/3}n^{2/3}}+O(n^{-1})\right], \label{leading-Asy-1}\\
\gamma_{n}=\frac{n^{-n-\frac \alpha 2-\beta-\frac 12}e^n}{\sqrt{2\pi}}
\left[1 +\frac{\sigma(s)}{2^{1/3}n^{1/3}} +\frac{\sigma^2(s)+2(\alpha+2\beta+1)u(s) }{2^{5/3}n^{2/3}}+O(n^{-1})\right],\label{leading-Asy-2}\\
\frac{a_n}{n}=2-\frac {2^{4/3}u(s)}{n^{2/3}}+O\left(n^{-1}\right),\label{an-asymptotic}\\
\frac{b_n}{n}=1-\frac {2^{1/3}u(s)} {n^{2/3}} +O\left (n^{-1}\right ),\label{bn-asymptotic}
\end{gather}
where  $\sigma(s)=\sigma(s;\beta,\omega,\mu)$  is the analytic solutions to  the   Jimbo-Miwa-Okamoto $\sigma$-form of the Painlev\'{e} II equation \eqref{sigma form} and   $u(s)=u(s;\beta,\omega,\mu)$ is the analytic solutions to  the $\mathrm{P}_{34}$ equation \eqref{pXXXV}.
\end{thm}

 When $\beta=0$, the asymptotic formulas of the recurrence coefficients have been studied by Lyu and Chen; see \cite{LyuChen}. It is seen that the  parameter $\omega$ in \eqref{weight} is not appeared in  equations \eqref{sigma form} and \eqref{pXXXV}. A similar phenomenon can be found in \cite{BasorChen,XuZhao,LyuChenFan}.

\begin{thm}\label{Asymptotic of monic}
Let $\pi_n(z)$ be the monic orthogonal polynomials defined as  \eqref{gamma}, then
\begin{align}\label{monic}
&\pi_n(4nz)=\frac{n^ne^{2n\left(z-\sqrt{z(z-1)}\right)}\left(\sqrt{z}+\sqrt{z-1}\right)^{2n+\alpha+2\beta}}{2^{\alpha+2\beta}e^n}
\left\{\frac{\sigma(s)}{2^{\frac 43}z^{\frac 14+\frac \alpha 2}(z-1)^{\frac 34+\beta}}\frac {1}{n^{\frac 13}}\right.\nonumber\\
&\left. +\left[\frac{2\sigma^2(s)-2u(s)-s}{2^{\frac{11}{3}}(\sqrt{z}+\sqrt{z-1})z^{\frac 14+\frac \alpha2}(z-1)^{\frac 54+\beta}}+
\frac{(\alpha+2\beta)(2u(s)+s)}{2^{\frac{8}{3}}z^{\frac 14+\frac \alpha2}(z-1)^{\frac 34+\beta}}\right]\frac1{n^{\frac 23}}+O\left(\frac 1n\right)\right\}
\end{align}
is valid uniformly on compact subsets of $\mathbb{C}\backslash [0,1]$, where $\arg z\in(-\pi,\pi)$ and $\arg(z-1)\in(-\pi,\pi)$.
\end{thm}

The paper is organized as follows. In Section \ref{sec:2}, we state a nonlinear steepest descent analysis of the RH problem for the orthogonal polynomial associated with the weight \eqref{weight}. In Section \ref{sec:4}, we give the proofs of the theorems \ref{Asymptotic of leading coefficients} and \ref{Asymptotic of monic}.
\section{Nonlinear steepest descent  analysis of orthogonal polynomials} \label{sec:2}
In this section, we take a nonlinear steepest descent analysis of the Rieman-Hilbert problem for $Y(z;\alpha,\beta,\omega,\mu)$ associated with the weight \eqref{weight} via the Deift-Zhou method; see \cite{DeiftZhouA},  \cite{DeiftZhouU}, and \cite{DeiftZhouS1999}.
\subsection{Riemann-Hilbert problem for orthogonal polynomials }
\begin{description}
  \item(Y1)~~  $Y(z;\alpha,\beta,\omega,\mu)$($Y(z)$  for short) is analytic in
  $\mathbb{C}\backslash [0,+\infty)$.

  \item(Y2)~~  $Y(z)$  satisfies the following jump condition
  \begin{equation}\label{}
  Y_+(x)=Y_-(x) \left(
                               \begin{array}{cc}
                                 1 & w(x) \\
                                 0 & 1 \\
                                 \end{array}
                             \right),
\quad x\in (0,+\infty),\end{equation} where the weight $w(x)$ is given in \eqref{weight}.

  \item(Y3)~~  The asymptotic behavior of $Y(z)$ at infinity is
  \begin{equation}\label{Y-infinity}Y(z)=\left (I+\frac {Y_{1}}{z}+\frac {Y_{2}}{z^2}+O(z^{-3})\right )\left(
                               \begin{array}{cc}
                                 z^n & 0 \\
                                 0 & z^{-n} \\
                               \end{array}
                             \right),\quad  z\rightarrow
                             \infty.\end{equation}
\item(Y4)~~ The behavior of $Y(z)$ at the point $0$ is
\begin{equation}\label{}Y(z)= \left\{\begin{array}{ll}
O\left(
                               \begin{array}{cc}
                                 1 & |z|^{\alpha} \\
                                 1 & |z|^{\alpha} \\
                               \end{array}\right), & z\rightarrow0,\quad \alpha <0,\\
O\left(
                               \begin{array}{cc}
                                 1 & \log|z| \\
                                 1 & \log|z| \\
                               \end{array}\right), & z\rightarrow0,\quad \alpha =0,\\
O\left(
                               \begin{array}{cc}
                                 1 & 1 \\
                                 1 & 1 \\
                               \end{array}\right), & z\rightarrow0,\quad \alpha >0.
 \end{array}\right.\end{equation}
\item(Y5)~~ The behavior of $Y(z)$ at the point $\mu$ is
\begin{equation}\label{}Y(z)= \left\{\begin{array}{ll}
O\left(
                               \begin{array}{cc}
                                 1 & |z-\mu|^{2\beta} \\
                                 1 & |z-\mu|^{2\beta} \\
                               \end{array}\right), & z\rightarrow\mu,\quad -\frac 12<\beta <0,\\
O\left(
                               \begin{array}{cc}
                                 1 & \log|z-\mu| \\
                                 1 & \log|z-\mu| \\
                               \end{array}\right), & z\rightarrow\mu,\quad \beta =0,\\
O\left(
                               \begin{array}{cc}
                                 1 & 1 \\
                                 1 & 1 \\
                               \end{array}\right), & z\rightarrow\mu,\quad \beta >0.
 \end{array}\right.\end{equation}
\end{description}

According to \cite{Fokas}, the unique solution of the Riemann-Hilbert problem for $Y(z)$ is given by
\begin{equation}\label{Y}
Y(z)= \left (\begin{array}{cc}
\pi_n(z)& \frac 1{2\pi i}\int_0^{+\infty} \frac{\pi_n(x) w(x)}{x-z}dx\\
-2\pi i \gamma_{n-1}^2 \;\pi_{n-1}(z)& - \gamma_{n-1}^2\;
\int_0^{+\infty} \frac{\pi_{n-1}(x) w(x)}{x-z}dx\end{array} \right ),
\end{equation}
where  the monic orthogonal polynomial $\pi_n(z)$ and leading coefficients $\gamma_n$ are defined as \eqref{gamma} and \eqref{recurrence relations} respectively.

\subsection{The  first transformation: $Y\rightarrow T$}
Before making the transformation, there is need to introduce some auxiliary functions. The density function respond to the external field $4x$ with $x>0$ is
$\varphi(x)=\frac 2 {\pi } \sqrt{\frac{1-x}{x}},~x\in(0,1),$
see \cite[(3.14)]{QiuWong} and \cite{Vanlessen,XuDanZhao}. Then we define
\begin{gather}g(z)=\int_{0}^{1}  \ln(z-x) \varphi(x) dx, \label{g}\\
\phi(z)=2\int_{1}^z\sqrt{\frac{1-x}{x}}dx=2\left[\sqrt{z(z-1)}-\ln\left(\sqrt{z-1}+\sqrt{z}\right)\right]
\label{phi}\end{gather}
with $\arg z\in (-\pi, \pi)$, $\arg(z-1)\in (-\pi, \pi)$ and $z\in\mathbb{C}\backslash(-\infty,1]$. The $g(z)$ and $\phi(z)$ are related by the  condition
$2g(z)+2\phi(z)-4z-l=0, ~ z\in \mathbb{C}\backslash(-\infty,1]$,
where $l=-2(1+2\ln 2)$; see \cite{QiuWong} for more details.

Now, we introduce the first transformation
 \begin{equation}\label{Y-T}T(z)= (4n)^{-(n+\frac\alpha 2+\beta)\sigma_3}e^{-\frac 1 2 nl \sigma_3} Y(4nz) e^ {n \left (\frac 1 2 l
-g(z)\right )\sigma_3}(4n)^{(\frac {\alpha}2+\beta) \sigma_3},~ z\in
\mathbb{C}\backslash [0,+\infty), \end{equation}
 where the Pauli matrices $\sigma_1=\begin{pmatrix}\begin{smallmatrix} 0 & 1\\ 1 & 0\end{smallmatrix}\end{pmatrix}$, $\sigma_2=\begin{pmatrix}\begin{smallmatrix} 0 & -i\\ i & 0\end{smallmatrix}\end{pmatrix}$, $\sigma_3 =\begin{pmatrix}\begin{smallmatrix} 1 & 0\\ 0 & -1\end{smallmatrix}\end{pmatrix}$. Then  $T(z)$ is analytic in $\mathbb{C}\backslash [0,+\infty)$,
and the behavior of $T(z)$ at infinity is $T(z)=I+O(1/z),~z\rightarrow\infty.$
\subsection{The  second transformation: $T\rightarrow S$}

We introduce  the transformation as
\begin{equation}\label{T-S}
S(z)=\left \{
\begin{array}{ll}
  T(z), & \mbox{for $z$ outside the lens region,}
  \\ &\\
  T(z) \left( \begin{array}{cc}
                                 1 & 0 \\
                                  -z^{-\alpha}(t-z)^{-2\beta} e^{2n\phi(z)} & 1 \\
                               \end{array}
                             \right) , & \mbox{ for $z$ in the upper lens region,}\\ &\\
T(z) \left( \begin{array}{cc}
                                 1 & 0 \\
                                  z^{-\alpha}(t-z)^{-2\beta} e^{2n\phi(z)} & 1 \\
                               \end{array}
                             \right) , & \mbox{ for $z$ in the lower lens region. }
\end{array}\right .\end{equation}
\begin{figure}[h]
\begin{center}
\includegraphics[scale=0.4,bb=0 0 554 156]{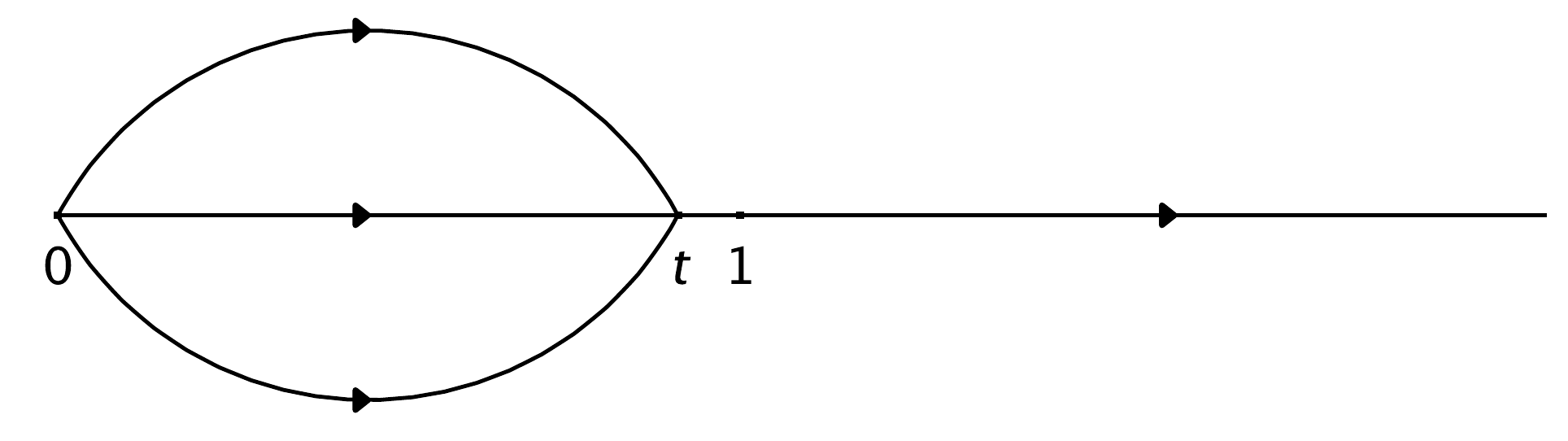}
\caption{The jump contours $\Sigma_{S}$ for the RH problem $S$. }
\label{contour-S} \end{center}
\end{figure}
Then
\begin{description}
  \item(S1)~~  $S(z)$ is analytic in
  $\mathbb{C}\backslash \Sigma_{S}$, where the contours $\Sigma_{S} $ are depicted in Figure \ref{contour-S}.

  \item(S2)~~  $S(z)$ satisfies the following jump conditions
  \begin{equation}\label{Jump-S}S_+(z)=S_-(z)
\begin{array}{ll}
  \left(
     \begin{array}{cc}
                                 1 & 0 \\
                                  z^{-\alpha}(t-z)^{-2\beta} e^{2n\phi(z)} & 1 \\
                               \end{array}
                             \right),&   z~\mbox{on lens},
\end{array}
\end{equation}
  and for $t<1$
\begin{equation}\label{Jump-S-1}S_+(x)=S_-(x)\left\{
\begin{array}{ll}
  \left(
                               \begin{array}{cc}
                                 0 & x^\alpha|x-t|^{2\beta} \\
                                 -x^{-\alpha}|x-t|^{-2\beta} & 0 \\
                               \end{array}
                             \right),&   x\in(0, t ), \\[.4cm]

\left(
 \begin{array}{cc}
                                e^{2n\phi_+(x)} & \omega x^\alpha|x-t|^{2\beta} \\
                                 0 & e^{2n\phi_-(x)}\\
                               \end{array}
                             \right),&   x\in(t,1),\\[.4cm]

\left(
                               \begin{array}{cc}
                                 1 & \omega x^\alpha|x-t|^{2\beta}e^{-2n\phi(x)} \\
                                 0 & 1 \\
                               \end{array}
                             \right),&   x\in(1,+\infty),
\end{array}
\right.
\end{equation}
and for $t>1$
\begin{equation}\label{Jump-S-2}S_+(x)=S_-(x)\left\{
\begin{array}{ll}
  \left(
                               \begin{array}{cc}
                                 0 & x^\alpha|x-t|^{2\beta} \\
                                 -x^{-\alpha}|x-t|^{-2\beta} & 0 \\
                               \end{array}
                             \right),&   x\in(0, 1),\\
  \left(
                               \begin{array}{cc}
                                 0 & x^\alpha|x-t|^{2\beta} e^{-2n\phi(x)}\\
                                 -x^{-\alpha}|x-t|^{-2\beta}e^{2n\phi(x)} & 0 \\
                               \end{array}
                             \right),&   x\in(1,t),\\
\left(
                               \begin{array}{cc}
                                 1 & \omega x^\alpha|x-t|^{2\beta}e^{-2n\phi(x)} \\
                                 0 & 1 \\
                               \end{array}
                             \right),&   x\in(t,+\infty).

\end{array}\right.
\end{equation}
  \item(S3)~~  The behavior of $S(z)$ at infinity is $S(z)=I+O(1/z),~z\rightarrow
                             \infty.$
\item(S4)~~ The behaviors of $S(z)$ at $z=0$ and $z=t$ are the same as those of $T(z)$ .
\end{description}

\subsection{Global Parametrix}
It is seen that the jump  matrix for $S$ tends to the identity on $\Sigma_S\backslash [0,t]$ as  $n\rightarrow\infty$,  the limiting Riemann-Hilbert problem $N(z)$ can be as follows:
\begin{description}
\item(N1)~~  $N(z)$ is analytic in  $\mathbb{C}\backslash
[0,t]$.
\item(N2)~~   $N(z)$ satisfies the jump condition
\begin{equation}\label{Jump-N} N_{+}(x)=N_{-}(x)\left(
       \begin{array}{cc}
       0 & x^\alpha|x-t|^{2\beta} \\
       -x^{-\alpha}|x-t|^{-2\beta} & 0 \\
       \end{array}
       \right),\quad x\in  (0,t).\end{equation}
\item(N3)~~    As $z\rightarrow\infty$, $ N(z)= I+O(z^{-1})$.
  \end{description}

It is seen that a solution of the RH problem for $N(z)$ can be constructed  as
\begin{equation}\label{N-expression}
 N(z) =d(\infty)^{\sigma_3} N_0(z) d(z)^{-\sigma_3},
 \end{equation}
 where
  $    N_0(z)=
\left(
  \begin{smallmatrix}
    \frac {\rho(z) + \rho^{-1}(z)} {2}&\frac {\rho(z) - \rho^{-1}(z)} {2i} \\
    -\frac {\rho(z) - \rho^{-1}(z)} {2i} &\frac {\rho(z) + \rho^{-1}(z)} {2} \\
  \end{smallmatrix}
\right), ~
 \rho(z)=\left ( \frac {z-t}{z} \right )^{1/ 4}, $
 and $d(z)$  is the Szeg\H{o} function,  which is analytic  in $\mathbb{C}\setminus [0,t]$ and satisfies
  $d_+(x)d_-(x)=x^\alpha|x-t|^{2\beta}, x\in(0,t)$; for Szeg\H{o} function see \cite{Szego1975}. Then $d(z)$ can be written  explicitly
  as
\begin{equation}\label{szego function}d(z)=z^{\frac \alpha 2}
(z-t)^\beta\left(\frac{2z-t+2\sqrt{z(z-t)}}{t}\right)^{-(\frac \alpha 2+\beta)}
\end{equation} with $d(\infty)=
\left (t/4\right )^{\alpha/2+\beta}$, where $\arg z\in (-\pi, \pi)$, $\arg(z-t)\in (-\pi, \pi)$ and the logarithm function $\ln z$ takes the principle branch such that $\arg z\in(-\pi,\pi)$; see \cite{KuijlaarsMcLaughlinAsscheVanlessen} for a relevant construction.


\subsection{Local parametrix at $z=t$}
The construction of a local parametrix $P^{(0)}(z)$ at a neighborhood  $U(0,r)$ can be refered to \cite{Vanlessen}, where the center of $U(0,r)$ is  $0$ with the radius $r$, here $r$ is small and positive. Hence, We concentrate on $z=t$.

Note that $ t= \frac {\mu}{4n}=1+\frac {s}{2^{2/3}n^{2/3}}$, we construct a local parametrix $P^{(1)}(z)$  at a neighborhood  $U(1,r)$, centered at $1$  with the radius $r$. $P^{(1)}(z)$  is analytic in $U(1,r) \backslash  \Sigma_{S}$, satisfies the same jump condition as that of $S(z)$ $\Sigma_{S}\cap U(1,r)$, $P^{(1)}(z)$ satisfies the matching condition
$P^{(1)}(z)N^{-1}(z)=I+ O\left (n^{-1/3}\right )$
on $\partial U(1,r)$, and has the same behavior at $z=t$ as that of $S(z)$.

In order to construct the $P^{(1)}(z)$, we need to introduce the
RH problem for  the $\mathrm{P}_{34}$ equation expressed as follows; see \cite{ItsKuijlaarsOstensson2008}, \cite{ItsKuijlaarsOstensson2009},
and \cite{WuXuZhao}. Let $\Phi(\zeta):=\Phi(\zeta;s)$, Then

\begin{figure}[h]
\begin{center}
\includegraphics[width=6cm,height=4cm,bb=0 0 555 555]{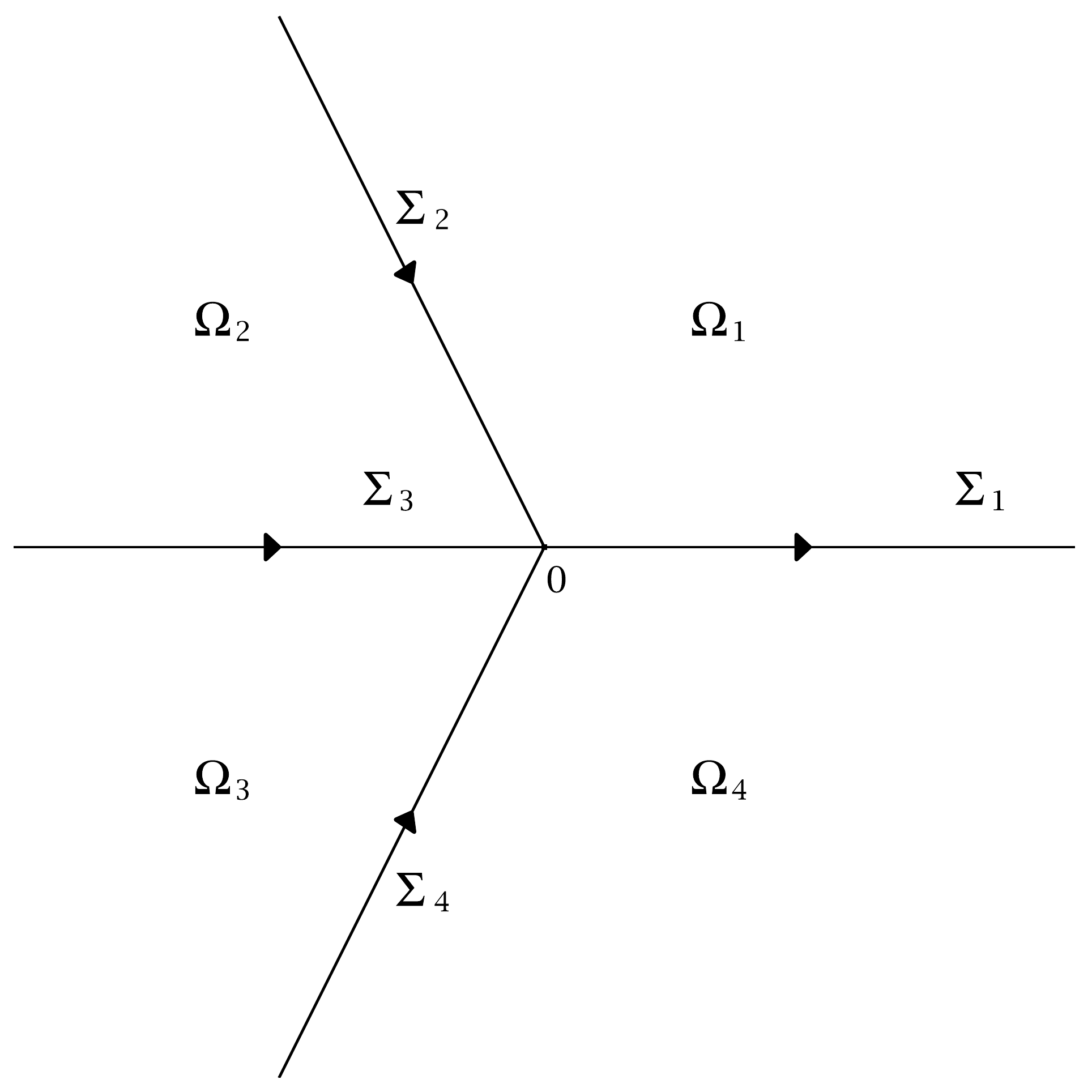}
\caption{The jump contours $\Sigma_i$ and regions for $\Phi$, with $\Sigma_1=\left\{\arg \zeta =0\right\}$, $\Sigma_2=\left\{\arg \zeta =\frac 23\pi\right\}$, $\Sigma_3=\left\{\arg \zeta =\pi\right\}$ and $\Sigma_4=\left\{\arg \zeta =-\frac 23\pi\right\}$.}
\label{fig1} \end{center}
\end{figure}

\begin{description}
  \item($\Phi1$)~~  $\Phi(\zeta)$  is analytic in
  $\mathbb{C}\backslash\cup_{i=1}^{4}\Sigma_i$; for the  contours $\Sigma_i$ cf. Figure \ref{fig1}.

  \item($\Phi2$)~~  $\Phi(\zeta)$  satisfies the  jump conditions
  \begin{gather}\label{Psi-jump}
  \Phi_+ (\zeta)=\Phi_- (\zeta)
 \left(
                               \begin{array}{cc}
                                 1 & \omega \\
                                 0 & 1 \\
                                 \end{array}
                             \right),
 \quad  \zeta \in {\Sigma}_1 ,\\
  \Phi_+(\zeta)=\Phi_-(\zeta)\left(
                               \begin{array}{cc}
                                 1 & 0 \\
                                 e^{2\beta\pi i} & 1 \\
                                 \end{array}
                             \right),  \quad   \zeta \in {\Sigma}_2 ,\\
\Phi_+ (\zeta)=\Phi_- (\zeta)
 \left(
                               \begin{array}{cc}
                                 0 & 1 \\
                                 -1 & 0 \\
                                 \end{array}
                             \right),
 \quad \  \zeta \in {\Sigma}_3 ,\\
  \Phi_+(\zeta)=\Phi_-(\zeta)\left(
                               \begin{array}{cc}
                                 1 & 0 \\
                                 e^{-2\beta\pi i} & 1 \\
                                 \end{array}
                             \right),  \quad   \zeta \in \Sigma_4.
  \end{gather}
\item($\Phi3$)~~     The behavior of $\Phi(z)$ at infinity is
\begin{align}\label{psi-infinity}
\Phi(\zeta)
&=\left(\begin{array}{cc} 1& 0\\ i\mathrm{m}_2(s) & 1 \end{array}\right)
     \left [I+\frac 1{\zeta}
     \left(\begin{array}{cc}\mathrm{m}_1(s) & -i\mathrm{m}_2(s)\\-i\mathrm{m}_3(s) & -\mathrm{m}_1(s)\end{array} \right)+O(\zeta^{-2})\right]\nonumber\\
  &\times \zeta^{-\frac{1}{4}\sigma_3}\frac{I+i\sigma_1}{\sqrt{2}} e^{-\theta(\zeta;s)\sigma_3} ,
 \end{align}
 where $\theta(\zeta;s)=\frac{2}{3}\zeta^{3/2}+s\zeta^{1/2}$, here $\arg \zeta\in(\pi,\pi)$.
\item($\Phi4$)~~ As $\zeta\rightarrow 0$, the behavior of $\Phi(z)$ is for $-\tfrac{1}{2}<\beta<0$
\begin{equation}\label{Psi0-1}
 \Phi(\zeta)=O\left(\begin{array}{cc} |\zeta|^\beta &  |\zeta|^\beta\\
 |\zeta|^\beta &  |\zeta|^\beta\end{array}\right)
 \end{equation}
 and for $\beta\geq0$
 \begin{equation}\label{Psi0-2}
\Phi(\zeta)=\left\{\begin{array}{ll}O\left(\begin{array}{cc} |\zeta|^\beta &  |\zeta|^{-\beta}\\
 |\zeta|^\beta &  |\zeta|^{-\beta}\end{array}\right), & \zeta\in \Omega_1\cup \Omega_4,\\
 O\left(\begin{array}{cc} |\zeta|^{-\beta} &  |\zeta|^{-\beta}\\
 |\zeta|^{-\beta} &  |\zeta|^{-\beta}\end{array}\right), & \zeta\in \Omega_2\cup \Omega_3.\end{array}\right.
\end{equation}
\end{description}

It is seen from \cite{ItsKuijlaarsOstensson2008,WuXuZhao} that for $\beta>-\frac 12$, $\omega\in \mathbb{C}\setminus (-\infty,0)$, there exits unique solution to the Riemann-Hilbert problem for $\Phi(\zeta;s)$, and
$\sigma(s)=\mathrm{m}_2(s)+\frac {s^2}{4}$ and
$u(s)=-\sigma'(s)=-\mathrm{m}_2'(s)-\frac s2$
satisfy the equations \eqref{sigma form} and \eqref{pXXXV} respectively.

Define  the conformal mapping as
\begin{equation}\label{f-t}
f(z)=\left (\frac{3}{2}\phi(z)\right )^{2/3}
\end{equation}
with $f(z)=2^{2/3}(z-1)+O(z-1)^2,~ z\rightarrow 1$, and
\begin{equation}\label{E}
E(z)=N(z)z^{\frac \alpha 2\sigma_3}(z-t)^{\beta \sigma_3}\frac{1}{\sqrt{2}}(I-i\sigma_1)
\left [n^{2/3}(f(z)-f(t))\right ]^{\sigma_3/4} .\end{equation}
Then, we introduce the following local parametix
\begin{equation}\label{parametrix}
P^{(1)}(z)=E(z)\Phi\left (n^{2/3}(f(z)-f(t));n^{2/3}f(t)\right )e^{n\phi(z)\sigma_3}(z-t)^{-\alpha \sigma_3}.\end{equation}
It is easily seen that the matching condition $P^{(1)}(z)N^{-1}(z)=I+ O\left (n^{-1/3}\right )$ is fulfilled.

\subsection{The final transformation: $S\rightarrow R$}
\begin{figure}[h]
\begin{center}
\includegraphics[scale=0.4,bb=1 0 554 149]{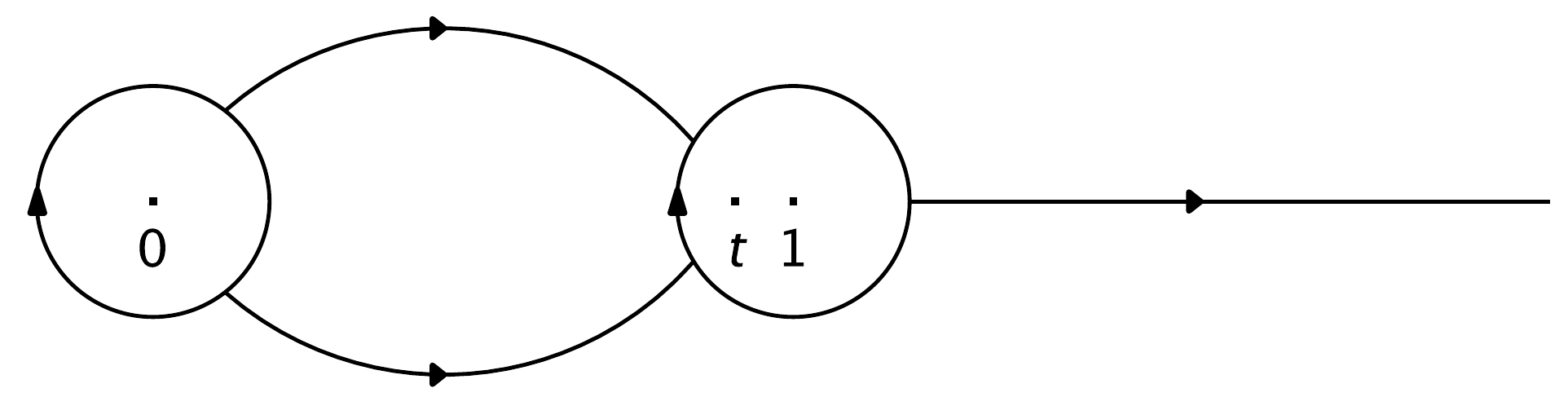}
\caption{The jump contours $\Sigma_R$ for the RH problem $R$. }
\label{contour-R} \end{center}
\end{figure}
The final transformation  is defined as follows
\begin{equation}\label{S-R}
R(z)=\left\{ \begin{array}{ll}
                S(z)N^{-1}(z), & z\in \mathbb{C}\backslash \left \{ U(0,r)\cup U(1,r)\cup \Sigma_S \right \},\\
               S(z) (P^{(0)})^{-1}(z), & z\in   U(0,r)\backslash \Sigma_{S},  \\
               S(z)  (P^{(1)})^{-1}(z), & z\in   U(1,r)\backslash
               \Sigma_{S} .
             \end{array}\right .
\end{equation}
Then, $R(z)$ is analytic in $\mathbb{C}\backslash \Sigma_R$, satisfies the jump conditions
$ R_+(z)=R_-(z)J_R(z), z\in \Sigma_R,$
and tends to $I$ as $z\rightarrow \infty$; for jump contours $\Sigma_R$ see Figure \ref{contour-R}. Note that as $n\rightarrow\infty$, $J_R(z)= 1+O(n^{-1/3})$ for $z \in \partial U(1,r)$.
Therefore, as $n\rightarrow\infty$, $R(z)=I+O(n^{-1/3})$,
where it is uniform for $z$ in whole complex plane.

\section{Proof of the main theorems}\label{sec:4}
In this section, we prove the theorem  \ref{Asymptotic of leading coefficients} and theorem  \ref{Asymptotic of leading coefficients}. A similar derivation can be found in \cite{XuZhao}, \cite{WuXuZhao} and \cite{KuijlaarsMcLaughlinAsscheVanlessen}.

\subsection {Proof of Theorem  \ref{Asymptotic of leading coefficients}}
It is seen, from \eqref{Y-infinity} and \eqref{Y}, that the leading coefficients of orthogonal polynomials are given by
 \begin{align}\label{leading-coefficients-R}
\gamma_{n-1}^2=-\frac{1}{2\pi i}(Y_1)_{21},\quad
\gamma_n^2=-\frac{1}{2\pi i}(Y_1)_{12}^{-1},
\end{align}
see \cite[(3.12)]{DeiftZhouS1999} and \cite[(5.7)]{BleherLiechty}. Tracing back the transformations \eqref{Y-T} \eqref{T-S} and \eqref{S-R} yields to
\begin{equation}\label{Y tracing back}
Y(4nz)=(4n)^{(n+\frac\alpha 2+\beta)\sigma_3}e^{\frac 12 nl\sigma_3}R(z)N(z)e^{ng(z)\sigma_3-\frac 12 nl \sigma_3}(4n)^{-(\frac {\alpha}2+\beta)\sigma_3}.
\end{equation}
Denote $N(z)$ and $R(z)$ as
\begin{gather}
N(z)=I+\frac{N_1}{z}+\frac {N_2}{z^2}+O(1/z^3),\quad z\rightarrow\infty,\label{N-expand}\\
R(z)=I+\frac{R_1}{z}+\frac {R_2}{z^2}+O(1/z^3),\quad z\rightarrow\infty.\label{R-expand}
\end{gather}
Note the fact that $e^{ng(z)\sigma_3}z^{-n\sigma_3}=I+G_2/z^2
+O(z^{-3})$ as $z\rightarrow\infty$, with $G_2=-\frac{n}{2}\int_{-1}^1 x^2\varphi(x)dx \sigma_3$,  we from \eqref{Y tracing back} have
\begin{gather}
Y_1=4n(4n)^{(n+\frac {\alpha}2+\beta)\sigma_3}e^{\frac 12 n l \sigma_3}(R_1+N_1)e^{-\frac 12 n l \sigma_3}(4n)^{-(n+\frac {\alpha}2+\beta)\sigma_3},\label{Y1}\\
Y_2=16n^2(4n)^{(n+\frac {\alpha}2+\beta)\sigma_3}e^{\frac 12 n l \sigma_3}(R_2+G_2+N_2+R_1N_1)e^{-\frac 12 n l \sigma_3}(4n)^{-(4n+\frac {\alpha}2+\beta)\sigma_3}.\label{Y2}\end{gather}

Next, we expand $N(z)$ and $R(z)$ as $z\rightarrow\infty$. From \eqref{N-expression}, we have
\begin{equation}\label{N-expand-1}
N(z)=I+\frac{N_1}{z}+\frac{N_2}{z^2}+O(z^{-3}),
\end{equation}
with
\begin{equation}\label{N-expand-2}
N_1=\left(\frac {t}{4}\right)^{(\frac\alpha 2+\beta)\sigma_3}\left[\frac { -t\sigma_2+(2\beta-\alpha)t\sigma_3}{4}\right]\left(\frac {t}{4}\right)^{-(\frac \alpha 2+\beta)\sigma_3},
\end{equation}
\begin{align}\label{N-expand-3}
N_2=\left(\frac {t}{4}\right)^{(\frac\alpha 2+\beta)\sigma_3}&\left[\frac {i (\alpha-2\beta)t^2\sigma_1}{16}-\frac { t^2\sigma_2}{8}+\frac { (10\beta-3\alpha)t^2 \sigma_3}{32}\right.\nonumber\\
&\left.\quad\quad+\frac {(\alpha^2-4\alpha\beta+4\beta^2+1)t^2I}{32}\right]\left(\frac {t}{4}\right)^{-(\frac \alpha 2+\beta)\sigma_3},
\end{align}
where the Pauli matrices $\sigma_1$, $\sigma_2$ and $\sigma_3$ are given in \eqref{Y-T}, and we have made use of the facts that as $z\rightarrow\infty$
\begin{gather*}\label{D-expand}
d(z)=\left(\frac {t}{4}\right)^{\frac \alpha 2+\beta}\left[1+\frac{(\alpha-2\beta)t}{4z}+\frac {(\alpha^2-4\alpha\beta+4\beta^2+3\alpha-10\beta)t^2)}{32z^2}+O(z^{-3})\right],\\
N_0(z)=I+\frac{-\frac{t}{4}\sigma_2}{z}+\frac{-\frac{t^2}{8}\sigma_2+ \frac{t^2}{32}I}{z^2}+O(z^{-3}).
\end{gather*}
Now we are in position to  estimate the behavior of $R(z)$ as $z\rightarrow\infty$. In view  of  \eqref{psi-infinity}, \eqref{E} and \eqref{parametrix}, 
we get the following  expression  of the jump for $R(z)$ on $\partial U(1, r)$ as $n\rightarrow\infty$
\begin{align}\label{R-jump-2}
J_R(z)= &N(z)z^{\frac {\alpha}{2}\sigma_3}(z-t)^{\beta \sigma_3} \nonumber\\
& \times\left[I+\frac{(F+\varphi_1)\sigma_3
-i\varphi_1\sigma_1}{n^{1/3}}+\frac{-\tau_3(z;t) \sigma_2+\tau_4(z;t) I
}{n^{2/3}}+O(n^{-1})\right]\nonumber\\
&\times(z-t)^{-\beta \sigma_3} z^{-\frac {\alpha}{2}\sigma_3}N^{-1}(z),
\end{align}
where
\begin{align}
&F(z;t,n)=\frac{s^2}{4(f(z)-f(t))^{1/2}},~
\varphi_1(z;t)=\frac{\mathrm{m}_2(s)}{2(f(z)-f(t))^{1/2}},\label{F-1}\\
&\varphi_2(z;t)=-\frac{\mathrm{m}^2(s)+\mathrm{m}_2'(s) }{2(f(z)-f(t))},~\tau_3(z;t)=F\varphi_1-\varphi_2,~
 \tau_4(z;t)=F\varphi_1+\frac {F^2} 2\label{r-3,4}.
\end{align}
Note the fact that $N_0(z)$ can be written in the following form
$N_0(z)=\frac{I-i\sigma_1}{\sqrt{2}} \rho(z)^{-\sigma_3}\frac{I+i\sigma_1}{\sqrt{2}},$
see \cite[(106)]{XuZhao}. We get on $\partial U(1, r)$
\begin{equation}\label{rort}
J_R(z) =d(\infty)^{\sigma_3}\frac{I-i\sigma_1}{\sqrt{2}} \left[I+\frac{\Theta_1(z)}{n^{1/3}}+\frac{\Theta_2(z)}{n^{2/3}}+O(n^{-1})\right]\frac{I+i\sigma_1}{\sqrt{2}} d(\infty)^{-\sigma_3},
\end{equation}
where
\begin{align*}
\Theta_1(z)
&=\tau_1(z;t)\sigma_++\tau_2(z;t)\sigma_-\nonumber\\
&~~~+i\varphi_1\left[\frac{i\xi^2(z;t) -i\xi^{-2}(z;t)}{2}\sigma_3-
\frac{\xi^2(z;t)+\xi^{-2}(z;t)-2}{2}\rho(z)^{-\sigma_3}\sigma_1\rho(z)^{\sigma_3}\right],\\
\Theta_2(z)
&=\tau_3(z;t) \sigma_3+\tau_4(z;t)I\\
& +\frac{\xi^2(z;t)+\xi^{-2}(z;t)-2}{2}\tau_3(z;t) \sigma_3+
\frac{i\xi^{2}(z;t)-i\xi^{-2}(z;t)}{2}\tau_3(z;t)\rho(z)^{-\sigma_3}\sigma_1\rho(z)^{\sigma_3},
\end{align*}
with $\sigma_+=\left(\begin{smallmatrix} 0 & 1\\ 0 & 0\end{smallmatrix}\right)$, $\sigma_-=\left(\begin{smallmatrix} 0 & 0\\ 1 & 0\end{smallmatrix}\right)$ and
\begin{gather}\label{r-1,2,xi}
\xi(z;t)=\frac {z^{\frac {\alpha}{2}}(z-t)^{\beta }}{d(z)},\quad
\tau_1(z;t)=-\rho^{-2}(z)(F+2\varphi_1)i,\quad \tau_2(z;t)=\rho^2(z) F i.
\end{gather}
We write
\begin{equation}\label{R asymptotic expansion}
R(z) =d(\infty)^{\sigma_3}\frac{I-i\sigma_1}{\sqrt{2}}
\left[I+\frac{R^{(1)}(z)}{n^{1/3}}+\frac{R^{(2)}(z)}{n^{2/3}}+O(n^{-1})\right]\frac{I+i\sigma_1}{\sqrt{2}}d(\infty)^{-\sigma_3}.
\end{equation}
Then,  $R^{(1)}(z)$ and   $R^{(2)}(z)$ satisfy the following scalar Riemann-Hilbert problems
\begin{gather}\label{}
R^{(1)}_+(z)-R^{(1)}_-(z)=\Theta_1(z),\quad
R^{(2)}_+(z)-R^{(2)}_-(z)=\Theta_2(z)+R^{(1)}_-(z)\Theta_1(z).
\end{gather}
Note that \eqref{r-1,2,xi} and \eqref{szego function}, one can obtain
\begin{gather*}\label{xi-expand}
   \xi(z;t)=1+\frac{(\alpha+2\beta)}{\sqrt{t}} \sqrt{z-t}+\frac{(\alpha+2\beta)^2}{2t}(z-t)+O((z-t)^{3/2}),\quad z\rightarrow t,\\
   \xi^2(z;t)+\xi^{-2}(z;t)-2=O(z-t),\quad z\rightarrow t\\   \xi^2(z;t)-\xi^{-2}(z;t)=\frac{4(\alpha+2\beta)}{\sqrt{t}}\sqrt{z-t}+O((z-t)^{3/2}),\quad z\rightarrow t.
\end{gather*}
It is seen, from \eqref{r-3,4} and \eqref{r-1,2,xi}, that $\tau_2(z;t)$  is analytic at $z=t$,  and $\tau_1(z;t)$, $\tau_3(z;t)$ and $ \tau_4(z;t)$ have a simple pole at $z=t$. Let $\mathrm{k}_2(s)=\lim\limits_{z\rightarrow t} \tau_2(z;t)$, $\mathrm{k}_j(s)=\Res\limits_{z=t} \tau_j(z;t)$ with $j=1,3,4.$ Hence, there have
\begin{gather}
 R^{(1)}(z)=\left\{\begin{array}{ll}
                   \frac{\mathrm{k}_1(s)}{z-t}\sigma_+-\Theta_1(z), &z\in U(1,r),\\[.3cm]
                    \frac{\mathrm{k}_1(s)}{z-t}\sigma_+, &z\not\in \overline{U(1,r)},
                  \end{array}
\right.\label{R-1}\\
 R^{(2)}(z)=\frac{\mathrm{k}_3(s)
\sigma_3+ \mathrm{k}_4(s) I+\mathrm{k}_5(s)\sigma_+-\mathrm{k}_1(s)\mathrm{k}_2(s)\sigma_-\sigma_+}{z-t} ,\ \ z\not\in
\overline{U(1,r)}.\label{R-2}\end{gather}
where $\mathrm{k}_5(s)$ is related to the  residue of the  $\xi(z;t)$ in $R^{(1)}_-(z)\Theta_1(z)$ and $\Theta_2(z)$. Employing \eqref{r-3,4} and \eqref{r-1,2,xi}, we have
\begin{equation} \label{k-i}
\begin{split}
 \mathrm{k}_1(s)&=-\frac{i\sigma(s)}{2^{1/3}} +O(n^{-\frac 23}),\quad
\mathrm{k}_2(s)=\frac{is^2} {2^{7/3}}+O(n^{-\frac 23}),\\
\mathrm{k}_3(s) &=\frac{4\sigma^2(s)-s^2\sigma(s)-4u(s)-2s}{2^{11/3}}
+O(n^{-\frac 23}),\\
\mathrm{k}_4(s) &=\frac{s^2\sigma(s)} {2^{11/3}}+O\left(n^{-\frac 23}\right),~
\mathrm{k}_5(s)=-2^{-\frac 23}i(\alpha+2\beta)\left(u(s)+\frac s2\right)+O\left(n^{-\frac 23}\right).
\end{split}\end{equation}
Expanding $\frac{1}{z-t}$ for $z$ at infinity, together with \eqref{R-1}, \eqref{R-2}, $R_j$ in \eqref{R-expand} and \eqref{R asymptotic expansion}, we get for $j=1,2$,
\begin{align}\label{Ri}
R_j&=d(\infty)^{\sigma_3}\left[\frac{\mathrm{k}_1(s)(i\sigma_3+\sigma_1)}{2n^{1/3}}+
\frac{\mathrm{k}_5(s)(i\sigma_3+\sigma_1)-\left(\mathrm{k}_1(s)\mathrm{k}_2(s)+2\mathrm{k}_3(s)\right )\sigma_2 }{2n^{2/3}} \right.\nonumber\\
&\left.\qquad\qquad\quad+O(n^{-1})  \right]d(\infty)^{-\sigma_3}.\end{align}
From $N_1$ in \eqref{N-expand-1},  $R_1$ in \eqref{Ri}, and $t=1+\frac{s}{2^{2/3}n^{2/3}}$, we have
 \begin{align}\label{leading-coefficients-R}
\gamma_{n-1}^2
&=-\frac{1}{2\pi i(4n)^{2n+\alpha+2\beta-1}e^{nl}}\left(\frac 4t\right)^{\alpha+2\beta}\nonumber\\
&\quad\times\left[-\frac{it}{4}+\frac{\mathrm{k}_1(s)}{2n^{1/3}}+
\frac{-i(2\mathrm{k}_3(s)+\mathrm{k}_1(s)\mathrm{k}_2(s))+\mathrm{k}_5(s)}{2n^{2/3}}+O(n^{-1})\right],
\end{align}
 \begin{align}\label{leading-coefficients-R}
\gamma_n^2
&=-\frac{1}{2\pi i(4n)^{2n+\alpha+2\beta-1}e^{nl}}\left(\frac 4t\right)^{\alpha+2\beta}\nonumber\\
&\quad\times\left[\frac{it}{4}+\frac{\mathrm{k}_1(s)}{2n^{1/3}}+
\frac{i(2\mathrm{k}_3(s)+\mathrm{k}_1(s)\mathrm{k}_2(s))+\mathrm{k}_5(s)}{2n^{2/3}}+O(n^{-1})\right].
\end{align}
Hence, from  \eqref{k-i}, we can complete the proof of \eqref{leading-Asy-1} and \eqref{leading-Asy-2} given in Theorem \ref{Asymptotic of leading coefficients}.

The recurrence coefficients $a_n$ as seen in \eqref{recurrence relations} can be expressed as
\begin{equation}\label{recurrence coefficients-a}
a_n=(Y_1)_{11}+\frac {(Y_2)_{12}}{(Y_1)_{12}},\quad b^2_n=(Y_1)_{12}(Y_1)_{21},
\end{equation}
see \cite[(3.12)]{DeiftZhouS1999}, where $Y_1$ and $Y_2$ are defined as \eqref{Y-infinity}. In view of \eqref{Y1} and \eqref{Y2}, there have
\begin{equation}\label{an-R-i}
\frac{a_n}{4n}=(N_1)_{11}+(R_1)_{11}+\frac{(N_2)_{12}+(R_2)_{12}+(R_1)_{11}(N_1)_{12}+(R_1)_{12}(N_1)_{22}}
{(N_1)_{12}+(R_1)_{12}}.
\end{equation}
In view of  \eqref{N-expand-1} and \eqref{Ri},  \eqref{k-i}  and $t=1+\frac{s}{2^{2/3}n^{2/3}}$, the equation \eqref{an-R-i} is simplified to
\begin{align}\label{an-asymptotic-proof}
\frac{a_n}{4n}&=\frac {t}{2}+\frac {\mathrm{k}_1(s)^2+\mathrm{k}_1(s)\mathrm{k}_2(s)+2\mathrm{k}_3(s)}{n^{2/3}}+O(n^{-1}).
\end{align}
Thus, \eqref{an-asymptotic} in theorem \ref{Asymptotic of leading coefficients} follows from \eqref{an-asymptotic-proof} and \eqref{k-i}. It follows from \eqref{Y1} that
\begin{equation}\label{b-n and R-i}
b_n^2=16n^2\left[(N_1)_{12}(N_1)_{21}+
(N_1)_{12}(R_1)_{21}+(N_1)_{21}(R_1)_{12}+
(R_1)_{12}(R_1)_{21} \right].
\end{equation}
From \eqref{N-expand-1}, \eqref{Ri}, \eqref{k-i}, and $t=1+\frac{s}{2^{2/3}n^{2/3}}$, we have
\begin{align}\label{bn-asymptotic-1}
b_{n} &  = 4n \displaystyle{\left[\frac {t^2}{16}+\frac { (2\mathrm{k}_3(s)+\mathrm{k}_1(s)\mathrm{k}_2(s))t+ \mathrm{k}_1 ^2(s) }{4n^{2/3}}
  +O(n^{-1})\right]^{1/2}}.
\end{align}
Note that \eqref{k-i}, the proof of \eqref{an-bn-asymptotic} in theorem \ref{Asymptotic of leading coefficients} is completed.
\subsection {Proof of Theorem \ref{Asymptotic of monic}}
From \eqref{Y tracing back} and \eqref{Y}, we have
\begin{equation}\label{}
\pi_n(4nz)=(4n)^ne^{ng(z)}(R_{11}N_{11}+R_{12}N_{21}).
\end{equation}
When $z\not\in \overline{U(1,r)}$, according to \eqref{R asymptotic expansion}, \eqref{R-1} and \eqref{R-2}, there have
\begin{align*}\label{R asymptotic expansion-1}
R(z)&=d(\infty)^{\sigma_3}\frac{I-i\sigma_1}{\sqrt{2}}
\left[I+\frac{\mathrm{k}_1(s)\sigma_+}{(z-t) n^{1/3}}\right.\nonumber\\
&\left.+\frac{(\mathrm{k}_3(s)\sigma_3+ \mathrm{k}_4(s)I+\mathrm{k}_5(s)\sigma_+-\mathrm{k}_1(s)\mathrm{k}_2(s)\sigma_-\sigma_+)}{(z-t)n^{2/3}}+O(n^{-1})\right]\frac{I+i\sigma_1}{\sqrt{2}}d(\infty)^{-\sigma_3}.
\end{align*}
Note that \eqref{N-expression} and \eqref{g}, we can complete the proof of  the theorem \ref{Asymptotic of monic}.
\section*{Acknowledgements}
This work was supported by the Science Foundation of Education Department of Jiangxi Province  under grant number GJJ170937 and the
National Natural Science Foundation of China under grant number 11801376.

\end{document}